\documentclass[12pt]{article}

\usepackage{amsmath,amsfonts,cd2}
\usepackage{amssymb}

\newtheorem{satz}{Theorem}[section]
\newtheorem{defi}[satz]{Definition}
\newtheorem{defis}[satz]{Definitions}
\newtheorem{kor}[satz]{Corollary}

\newtheorem{bems}[satz]{Remarks}
\newtheorem{prop}[satz]{Proposition}

\newtheorem{bspe}[satz]{Examples}

\newtheorem{rem}[satz]{Remark}

\newcounter{Roma}
\newcounter{Ara}
\newcounter{let}

\newenvironment{buch}{\begin{list}{(\alph{Ara})\hfill}{\usecounter{Ara}
\labelsep3mm \leftmargin1.1cm \labelwidth8mm}}{\end{list}}

\newenvironment{rom}{\begin{list}{(\roman{Roma})\hfill}{\usecounter{Roma}
\labelsep3mm \leftmargin1.5cm \labelwidth8mm}}{\end{list}}

\makeindex

\begin{document}

\newcommand{\nc}{\newcommand}

\nc{\mapco}{\,\colon\, }
\nc{\ab}{^{ab}}

\nc{\comment}[1]{
}

\nc{\catc}{{\C}}
\nc{\we}{\vee}

\nc{\hra}{\hookrightarrow}

\nc{\epi}{epimorphism}
\nc{\repi}{regular epimorphism}
\nc{\mono}{monomorphism}
\nc{\iso}{isomorphism}
\nc{\coker}[1]{\mbox{${\rm Coker}(#1)$}}
\nc{\Ker}[1]{\mbox{{\rm Ker}$(#1)$}}
\nc{\defgl}{\stackrel{def}{=}}

\nc{\V}{\vspace{3mm}}
\nc{\VV}{\vspace{4mm}}
\nc{\lra}{\longrightarrow}
\nc{\lla}{\longleftarrow}
\nc{\mr}[1]{ \stackrel{#1}{\lra} }
\nc{\ml}[1]{ \stackrel{#1}{\lla} }
\nc{\hmr}[1]{\hspace{2mm}\stackrel{#1}{\lra}\hspace{2mm}}
\nc{\hml}[1]{\hspace{2mm} \stackrel{#1}{\lla}\hspace{2mm}}
\nc{\N}{\noindent}
\nc{\st}{^{\prime}}
\nc{\ot}{\otimes}
\nc{\hcong}{ \hspace{2mm}\cong\hspace{2mm}  }
\nc{\hfbox}{\hfill$\Box$}
\nc{\REF}[1]{(\ref{#1})}

\def\Z{\ifmmode{Z\hskip -4.8pt Z} \else{\hbox{$Z\hskip -4.8pt Z$}}\fi}

\def\Q{\ifmmode{Q\hskip-5.0pt\vrule height6.0pt depth
0pt\hskip6pt}\else{\hbox{$Q\hskip-5.0pt\vrule height6.0pt depth
0pt\hskip6pt$}}\fi}

\newcommand{\NN}{\mbox{$I\!\!N$}}

\nc{\Ph}{\phantom{}}

\nc{\BE}{\begin{equation}}
\nc{\EE}{\end{equation}}

\nc{\dst}{\displaystyle}
\nc{\sst}{\scriptscriptstyle}
\nc{\ssst}{\scriptscriptstyle}
\nc{\proof}{\N{\bf Proof\,:}\quad}
\nc{\proofof}[1]{\N{\bf Proof of {#1}\,:}\quad}
\nc{\proofofthm}[1]{\N{\bf Proof of theorem \ref{#1}\,:}\quad}

\nc{\htt}[1]{^{\otimes #1}}

\newcommand{\ssur}[2]{\mbox{$#1 \!\to\!\!\!\!\!\to\! #2$}}

\def\mapsel#1{\mbox{$\rule[-5mm]{0mm}{12mm} \searrow\rlap
{$\vcenter{\makebox[0mm][r]{$\scriptstyle#1$\hspace{12.2mm}}}$}$}}

\def\surltop#1{\makebox[1cm]{\mbox{$\stackrel{#1}{\makebox[0mm]{$\lla$}\hspace{0.7mm}
\makebox[0mm]{$\lla$}}$}}}

\nc{\Sur}[1]{\mbox{$\:\stackrel{#1}{\lra\!\!\!\!\!\to\,}\:$}}

\def\INJ{\mbox{\mathsurround=0pt
\makebox[0mm][r]{\parbox{0mm}{\rule[-0.65mm]{0mm}{0.2mm}$\scriptscriptstyle>$}}
\makebox[0.7cm][l]{\parbox{0.7cm} {$\lra$}}}}

\def\Inj#1{\mbox{$\:\stackrel{#1}{\INJ}\:$}}

\nc{\Injup}[1]{\mapup{#1}}
\nc{\injup}[1]{\mapup{#1}}
\nc{\injdown}[1]{\mapdown{#1}}


\def\mapup#1{\mbox{$\rule[-1mm]{0cm}{0.7cm}  
\makebox[0mm][r]{\raisebox{0.2mm}{$\scriptstyle\phantom{\cong}$}\hspace{0.6mm}}
\bigg\uparrow\rlap{$\vcenter{\hbox{$\scriptstyle#1$}}$}$}}

\def\isoup#1{\mbox{$\rule[-1mm]{0cm}{0.7cm}  
\makebox[0mm][r]{\raisebox{0.2mm}{$\scriptstyle\cong$}\hspace{0.6mm}}
\bigg\uparrow\rlap{$\vcenter{\hbox{$\scriptstyle#1$}}$}$}}

\def\mapdown#1{\mbox{$\rule[-1mm]{0cm}{0.7cm}  
\makebox[0mm][r]{\raisebox{0.2mm}{$\scriptstyle $}\hspace{0.6mm}}
\bigg\downarrow\rlap{$\vcenter{\hbox{$\scriptstyle#1$}}$}$}}

\def\isodown#1{\mbox{$\rule[-1mm]{0cm}{0.7cm}  
\makebox[0mm][r]{\raisebox{0.2mm}{$\scriptstyle\cong$}\hspace{0.6mm}}
\bigg\downarrow\rlap{$\vcenter{\hbox{$\scriptstyle#1$}}$}$}}

\def\isor#1{\mbox{$\smash{\mathop{\longrightarrow}\limits^{\cong}_{#1}}$}}

\def\isol#1{\mbox{$\smash{\mathop{\longleftarrow}\limits^{\cong}_{#1}}$}}

\def\surdown#1{\makebox[0mm]{$\mapdown{\raisebox{0.4mm}{$\scriptstyle#1$}}$}
\makebox[0mm]{\raisebox{-2.15mm}{$\downarrow$}}}

\def\surup#1{\makebox[0mm]{$\mapup{\raisebox{0.4mm}{$\scriptstyle#1$}}$}
\makebox[0mm]{\raisebox{0.7mm}{$\mapup{}$}}}

\newcommand{\surr}[2]{@\cdhgeneric>->->\twoheadrightarrow>#1>#2>}

\newcommand{\surl}[2]{@\cdhgeneric>\twoheadleftarrow>->->#1>#2>}

\newcommand{\brokenr}[2]{
@\cdhgeneric>\raise2.8pt\hbox to3.5pt{\hrulefill}\mkern9mu>
\raise2.8pt\hbox to3.5pt{\hrulefill}\hbox to5pt{}>
\mkern-7.5mu\dashrightarrow>#1 >
#2>}

\newcommand{\brokenup}[1]{
@\cdvgeneric>\hat\cdot
>\raisebox{3pt}{$\vdots$}>
\vbox{\kern3.5pt\hbox{$\cdot$}\kern-3.5pt}
> >\hspace{1mm}#1>}

\newcommand{\functlr}[2]{
\raisebox{0.4pt}{$\hss\begin{CD}
@>\vbox{\hbox to 0pt{$\hss\begin{CD}@<#1<<\end{CD}\hss$}\vskip-2pt}
>#2 >
\end{CD}\hss$}
}

\newcommand{\maprr}[2]{
\raisebox{-0.9pt}{$\hss\begin{CD}
@>\vbox{\hbox to 0pt{$\hss\begin{CD}@>#1>>\end{CD}\hss$}\vskip-3pt}
>#2 >
\end{CD}\hss$}
}

\newcommand{\mapdd}[2]{
@\cdvstandard>\downarrow\hspace{1.5pt}\kern-4pt\hspace{1.5pt}\Big\downarrow>#1>#2>}

\newcommand{\mapud}[2]{
@\cdvstandard>\uparrow\kern-3pt\Big\downarrow>#1>#2>}

\newcommand{\sepi}[3]{\,\mbox{$#1\,: \ssur{#2}{#3}$}\,}

\nc{\auf}{\twoheadrightarrow}
 
\nc{\ruled}{\rule[-4mm]{0mm}{0mm}}

\nc{\qu}{quadratic}
\nc{\GBoGB}{G/BG\st \otimes G/BG\st}
\nc{\lstar}{_{\raisebox{-1mm}{$*$}}}
\nc{\sm}{\:{ \wedge}\:}

\nc{\Rmod}{${\bf R}$-module}

\nc{\map}[3]{\mbox{$#1 \mapco #2 \to #3$}}

\nc{\rond}{{\,\sst \circ\,}}

\nc{\ruleu}{\rule{0mm}{7mm}}

\nc{\T}[1]{\tilde{#1}}

\nc{\Imm}[1]{\mbox{${\rm Im}(#1)$}}

\nc{\tw}{\end{document}}

\nc{\QJG}{\frac{\dst I(G) J}{\dst I^2(G) J}}
\nc{\QJH}{\frac{\dst I(H) J}{\dst I^2(H) J}}
\nc{\otz}{\ot}
\nc{\UL}[2]{{\rm U}_{#1}{\rm L}(#2)}

\nc{\IRN}[1]{I_{R,\cal G}^{#1}(G)}

\nc{\ULG}[1]{{\rm U}_{#1}{\rm L}^{\cal G}(G)}
\nc{\UGH}[1]{{\rm U}_{#1}^{\cal GH}(G,H)}


\nc{\calG}{{\cal G}}
\nc{\AB}{^{AB}}

\nc{\IZ}{I}

 \nc{\UG}{\mbox{${\rm UL}^{\cal G}(G)$}}

 \nc{\UN}{\mbox{${\rm UL}^{\cal N}(N)$}}

 \nc{\UT}{\mbox{${\rm UL}^{\gamma}(T)$}}

 \nc{\UGn}[1]{\mbox{${\rm U}_{#1}{\rm L}^{\cal G}(G)$}}

 \nc{\UNn}[1]{\mbox{${\rm U}_{#1}{\rm L}^{\cal N}(N)$}}

 \nc{\UTn}[1]{\mbox{${\rm U}_{#1}{\rm L}^{\gamma}(T)$}}

\nc{\LG}{\mbox{$ {\rm L}^{\sst\cal G}(G)$}}
\nc{\LH}{\mbox{$ {\rm L}^{\sst\cal H}(H)$}}
\nc{\LN}{\mbox{$ {\rm L}^{\sst\cal N}(N)$}}

\nc{\LnG}[1]{\mbox{$ {\rm L}_{#1}^{\sst\cal G}(G)$}}

\nc{\LNn}[1]{\mbox{$ {\rm L}_{#1}^{\sst\cal N}(N)$}}
\nc{\LnN}[1]{\mbox{$ {\rm L}_{#1}^{\sst\cal N}(N)$}}

\nc{\KnD}[1]{{\cal K}_{#1}^{\Delta}}

\nc{\KnL}[1]{{\cal K}_{#1}^{\Lambda}}

\nc{\KnINL}[1]{{\cal K}_{#1}^{I(N)\Lambda}}




 \nc{\QnG}[1]{Q_{#1}(G)}

 \nc{\QnGT}[1]{Q_{#1}(G,T)}

 \nc{\QnGN}[1]{Q_{#1}(G,N)}

\nc{\Tor}[2]{{\rm Tor}_1^{Z\!\!\!Z}(#1\,,#2)}


\def\isor#1{\mbox{$\smash{\mathop{\longrightarrow}\limits^{\cong}_{#1}}$}}
\def\isol#1{\mbox{$\smash{\mathop{\longleftarrow}\limits^{\cong}_{#1}}$}}

\nc{\caln}{{\cal N}}
\nc{\UNN}[1]{{\rm U}^{\cal N}_{#1}(N,N)}
\nc{\ULN}[1]{{\rm U}_{#1}{\rm L}^{\cal N}(N)}


\nc{\calH}{{\cal H}}
\nc{\calK}{{\cal K}}

\nc{\ULH}[1]{{\rm U}_{#1}{\rm L}^{\cal H}(H)}


\nc{\hoplus}{\hspace{3mm}\oplus\hspace{3mm}}

\newcommand{\C}{\mbox{$ C\!\!\!\!I$}}
\newcommand{\RR}{\mbox{$I\!\!R$}}



\nc{\maprtop}[1]{\mr{#1}}
\nc{\mapltop}[1]{\ml{#1}}

\def\mapse#1{\mbox{$\Bigg{.\searrow}\rlap{$\vcenter{\hbox{$\scriptstyle#1$}}$}$}}

\def\mapsw#1{\mbox{$\Bigg{.\swarrow}\rlap{$\vcenter{\hbox{$\scriptstyle#1$}}$}$}}

\def\mapne#1{\mbox{$\Bigg{.\nearrow}\rlap{$\vcenter{\hbox{$\scriptstyle#1$}}$}$}}

\def\mapnw#1{\mbox{$\Bigg{.\nwarrow}\rlap{$\vcenter{\hbox{$\scriptstyle#1$}}$}$}}

\title{Quadratic maps between groups}

\author{{\normalsize\it Dedicated to Mamuka Jibladze on the occasion of his fiftieth birthday}\\ \\ Manfred Hartl}
\date{}
\maketitle


\begin{center}

\N LAMAV and FR CNRS 2956, 
ISTV2, \\ Universit\'e de Valenciennes et du Hainaut-Cambr\'esis, \\
Le Mont Houy,  59313 Valenciennes Cedex 9, France.\\ 
Email: Manfred.Hartl@univ-valenciennes.fr

\end{center}
\vspace{8mm}

\begin{abstract} The notion of \qu\ maps between arbitrary groups appeared at several places in the literature on \qu\ algebra. Here a unified extensive treatment of their properties is given; the relation with a relative version of Passi's polynomial maps and groups of degree 2 is established and used to study the structure of the latter.

\end{abstract}\vspace{5mm}

 \N{\large \bf Introduction.}\quad Polynomial maps appear in nilpotent group theory for a long time, originally in the form of rational (numerical) functions, for example in the Hall-Petrescu formula or the group law of torsionfree nilpotent groups when written with respect to a Mal'cev basis. An {\em intrinsic}\/ notion of polynomial maps from groups to {\em abelian}\/ groups was introduced by Passi \cite{Pa68}, together with a universal example $G\to P_n(G)$ where the abelian group $P_n(G)$ is called ``polynomial group". Passi's motivation came from the study of dimension subgroups; since then, his construction turned out to provide a key tool in the study of many other problems: in the theory of group schemes \cite{Br} as well as in the theory of nilpotent groups, concerning  their second (co)homology \cite{PolProp}, \cite{GoG},   automorphism groups or simplicial objects \cite{Diss}, \cite{Habil}. However, a need to study polynomial maps between {\em arbitrary}\/ groups comes from unstable homotopy theory; after Baues' \cite{Bametastable} and the author's \cite{Diplom} study of metastable homotopy  groups    and Moore spaces \cite{Moore}  the foundations of ``quadratic algebra" were layed  in \cite{BHP} where  a notion of quadratic maps  with non-abelian target group first appeared. Since then, in the steadily growing literature on quadratic algebra and its applications, various variants and properties  were exhibited when needed, in work of Baues, Jibladze, Muro, Pirashvili and the author  (most of these articles can be found on ArXiv). So the purpose of this paper is to provide a thorough unified treatment of what is called {\em weakly}\/ \qu\ maps, because of their good functorial behaviour; for brevity we drop the word ``weakly" in this paper. So several of our formulas and properties appear also elsewhere in the literature, but we include and prove most of them here for the sake of a coherent exposition. However, we here work in the slightly more general framework of quadratic maps \textit{relative to a subgroup}, inspired by Passi's study of relative dimension subgroups, see \cite{Pa}, lateron extended by Kuz'min \cite{Ku} and in \cite{HMP}. The relative viewpoint here leads to the construction of various categories of relative quadratic maps generalizing the ``quadratic envelope of the category of 2-step nilpotent groups"  introduced by Jibladze and Pirashvili \cite{Niq}. In particular, the category
  ${\bf CP}$ of central \qu\ pair maps introduced  here turns out to play a fundamental role in \qu\ algebra since it allows to refound the basic notions and to define modules over square ringoids instead of only square groups as in \cite{BHP}; this is work in progress and will be presented in \cite{ModQuad} and in a forthcoming book on \qu\ algebra jointly written with H.\ Gaudier, F.\ Goichot, B.\ Loiseau and T.\ Pirashvili.

 
 In section 1 we introduce and study relative \qu\ maps and their universal examples $G\to Q(G,B)$ which constitute a nonabelian version of the relative version of Passi's  polynomial maps and groups in degree 2. The latter are studied in section 2, first for arbitrary degree before focussing on the degree $2$ where the relative Passi groups $P_2(G,B)$
  turn out to be precisely the abelianization of the groups $Q(G,B)$; this fact allows to deduce their structure from the (easy) non-abelian case in a rather simple way, and to deduce several exact sequences for $P_2(G,B)$.
 
 Finally we note that most of the theory of this paper can be generalized to the abstract setting of semi-abelian categories \cite{QuadFunct}; in particular, all algebraic theories containing a group law as part of the structure (for example, algebras over any non-unitary $k$-linear  operad, like Lie algebras), admit a theory of quadratic maps. Also, a theory of \qu\ maps between modules is inaugurated in \cite{P2RM}, and a notion  of polynomial maps between non-abelian groups of arbitrary degree is introduced by the author, from an inductive viewpoint generalizing the one adopted in this paper; this is work in progress. There is, however, a different approach due to
 Leibman which also has interesting applications \cite{Leib1}, \cite{Leib2}; the precise relation between the two approaches  remains to be clarified. \V

\nc{\rad}[1]{\mbox{rad$(#1)$}}

\section{Quadratic maps between groups}

Throughout this paper, the symbols $G$ and $H$ denote groups.
The commutator of elements $a,b\in G$ is defined as 
$ [a,b] =aba^{-1}b^{-1}$,
and the conjugation is denoted by ${}^ab =aba^{-1}$. We write $G\ab =G/G'$ and $ab\mapco G\auf G\ab$ for the natural projection. Recall that the lower central series of $G$ is defined by $\gamma_1(G)=G$ and $\gamma_{i+1}(G)=[G,\gamma_i(G)]$. We say that $G$ is {\em $n$-step nilpotent}\/ if  $\gamma_{n+1}(G) =\{1\}$.  If $G$ is $2$-step nilpotent, the commutator map $[-,-]\mapco G\times G \to G$ is wellknown to be bilinear.\V

Let \map{f}{G}{H} be some function between {\em arbitrary}\/ groups. We shall, however, write the group law of $H$ {\em
additively}\/ since in many applications   $H$ is abelian, and in those in \cite{BHP}, \cite{ModQuad} where $H$ is
genuinely nonabelian, it is written additively anyway to match the conventions in homotopy theory which originally motivated these developments.

Define the {\em deviation function} or {\em cross effect}\/ of $f$ to be the map
\BE\label{defdf}  \map{d_f}{G\times G}{H} \quad \mbox{by} \quad d_f(a,b) = f(ab) -f(b)-f(a)\:.\EE
Furthermore, let $I_f$ resp.\ $D_f$ denote the subgroup of $H$ generated by \Imm{f} resp.\ \Imm{d_f}. 

\begin{defis}\label{qulindef} \rm We say that $f$ as above is \V

  (a) {\em linear}\/ if $d_f=0$, i.e., $f$ is a group homomorphism;\V

  (b) {\em quadratic}\/ if $d_f$ is {\em bilinear}\/ and $D_f$ is {\em central}\/ in $I_f$, or more explicitely,
$\forall a,b,c \in G$, $[d_f(a,b),f(c)] = 0$.\V

  This definition of \qu\ maps first appears in \cite{Niq} under the same of weakly \qu\ maps. 
  
  Note that linear maps are quadratic. We denote by ${\rm Quad}(G,H)$ the set of quadratic maps from $G$ to $H$.\V

We also need a relative version of quadratic maps, as follows. Let $B$ be a subgroup of $G$. Then we say that  $f$ as above is {\em quadratic
relative B}\/ if $f$ is quadratic and $d_f(B\times G) = d_f(G\times B) = 0$. Note that $f$ is \qu\ iff it is \qu\ relative $\{1\}$. We define the {\em radical}\/ of a quadratic map $f$ to be the set ${\rm rad}(f)$ consisting of all elements $a$ of $G$ such that for all $b\in G$, $d_f(a,b)=d_f(b,a) =0$. Note that \rad{f} is a subgroup of $G$ and that $G' \subset \rad{f}$ since  $d_f$ is a bilinear map taking values in the abelian group $D_f$. In particular, \rad{f} is normal and $G/\rad{f}$ is abelian. It is also clear that \rad{f} is the largest subgroup $B$ of $G$ such that $f$ is quadratic relative $B$.

\end{defis}

In the following proposition we collect the basic properties of quadratic maps which are easily deduced from the definition.\V

\begin{prop}\label{qumapsprops} Let $f\mapco G \to H$ be a quadratic map relative some subgroup $B$ of $G$. 
\begin{buch}

\item  One has the following identities for $a,b\in G$:
\BE\label{f(ab)=}  f(ab) \hspace{2mm}=\hspace{2mm} d_f(a,b) + f(a) + f(b) \hspace{2mm}=\hspace{2mm} f(a) + f(b)
+ d_f(a,b)\EE 
\BE\label{df(a,b)=}  d_f(a,b) = {}-f(a) + f(ab) -f(b) \:;\EE  

\item $f$ is normalized, i.e., $f(1)=0$;

\item the restriction of $f$ to $B$ is linear, whence $f(B)$ is a subgroup of $H$;

\item there is a canonical linear map 
 \BE\label{wf} \map{w_f}{G/BG\st \otimes G/BG\st}{H} \EE
such that $w_f(\bar{a}\ot \bar{b} ) = d_f(a,b)$ for $a,b\in G$.\hfbox

\end{buch}
\end{prop}\V











\begin{bspe}\label{Exqumaps} \rm (0) Let $R$ be any ring and $a,b\in R$. Then the function $f\,\colon R\to R$, $f(x) =ax^2+bx$, is a quadratic map
between additive groups. More generally, quadratic forms on vector spaces are \qu\ maps.\V

(1) For any subgroup $B$ of $G$,  the 2-fold diagonal map $\delta^2\,\colon G\to \GBoGB$, $\delta(a) =
\bar{a} \ot
\bar{a}$, is \qu\ relative $B$.\V

(2) Let $L$ be a 2-step nilpotent Lie algebra  (i.e., $[[L,L],L] =0$) over   $\Z[\frac{1}{2}]$.    
 Then a multiplicative group law  on $L$ is defined by the truncated
Campbell-Baker-Hausdorff-formula, i.e., $x\rond y = x+y+\frac{1}{2}[x,y]$ for $x,y\in L$. In fact, $(L,\rond)$ is a   
uniquely $2$-divisible 2-step nilpotent group, and this construction provides a functorial equivalence between groups of this type
and 2-step nilpotent Lie algebras over $\Z[\frac{1}{2}]$. This is a special case of a more general result of
Lazard
\cite{La},
 providing ``abelian models" for sufficiently divisible nilpotent groups. Now in the 2-step nilpotent case above,
the identity map \map{id}{(L,\rond)}{(L,+)} is {\em \qu}\/ with $d_{id}(x,y) = 
\frac{1}{2}[x,y]$.

We note that generalizing the above equivalence of Lazard we constructed functorial abelian models for {\em arbitrary}\/  2-step nilpotent groups (actually, for central group extensions with abelian cokernel), cf.\ \cite{Diss} and also \cite{Habil}; this construction is based on the properties of relative quadratic maps with values in abelian groups, see section 2 below.\V


 (3) Let $k$ be some commutative ring with unit and $G$ be the subgroup $1+Tk[[T]]$ of the group of units of the algebra
$k[[T]]$ of power series over $k$. For $n\ge 0$, let $c_n\,\colon\,G \to k$, $c_n(\sum_{i\ge 0} a_iT^i) = a_n$. Then $c_1$ is
linear and $c_2$ is quadratic since $c_2(fg) = c_2(f) + c_2(g) + c_1(f)c_1(g)$, $f,g\in G$.\V

 (4) Let $\underline{\Sigma}(n,3n-3)$ denote the pointed homotopy category of suspensions $\Sigma X $ of topological spaces $X$ such that 
$\Sigma X $ is an $(n-1)$-connected $(3n-3)$-dimensional CW-space. Then for $\Sigma X , \Sigma Y \in \underline{\Sigma}(n,3n-3)$ the second James-Hopf-invariant 
  \[ \gamma_2\mapco [\Sigma X, \Sigma Y ] \hmr{} [\Sigma X, \Sigma Y  \sm Y] \]
 (see \cite{Wh}) is a \qu\ map where the group structure on  $[\Sigma X, Z]$ is induced by the cogroup structure of $ \Sigma X$. 
  Moreover, if also $\Sigma Z \in \underline{\Sigma}(n,3n-3)$ and $f\mapco \Sigma X \to \Sigma Y$ is a continous map then the map
    \[ [f]^* \mapco  [\Sigma Y, \Sigma Z ] \hmr{}  [\Sigma X, \Sigma Z ]  \]
    is a \qu\ map, see \cite[Appendix]{Bametastable}. This example gave rise to  the notion  of quadratic categories, see \cite{BHP}.

\end{bspe}

More examples appear in the following proposition showing that \qu\ maps are intimely related to 2-step nilpotent groups.

\begin{prop}\label{quadnil2} Let $G$ be a group and $n\ge 2$. Then the following properties are equivalent.\V

(1) $G$ is 2-step nilpotent.\V

(2) The $(n-1)$-fold multiplication map \map{\mu_{n-1}}{G^n}{G}, $\mu_{n-1}(a_1,\ldots,a_n)= a_1 \cdots a_n$, is a \qu\ map.\V

(3) For all groups $K$ and all linear maps \map{f_1,\ldots,f_n}{K}{G}, the product map \map{f_1 \cdots f_n}{K}{G}, $x\mapsto
f_1(x) \cdots f_n(x)$, is \qu.\V

(4) The map \map{2_G}{G}{G}, $a\mapsto a^2$, is \qu.\V

\N Each of these implies\V

(5) The map \map{n_G}{G}{G}, $a\mapsto a^n$, is \qu.\V
\end{prop}

Note that property (2) neatly generalizes the often useful fact that $G$ is abelian iff $\mu_1$ is linear.\V

\N{\bf Proof\,:}  We first note that if $G$ is 2-step nilpotent the following relations hold for $a=(a_1,\ldots,a_n)$ and
$b=(b_1,\ldots,b_n)$ in $G^n$:
\begin{eqnarray}
\label{nil2formel1} a_1b_1 \cdots a_nb_n &=& \prod_{1\le i < j \le n} [b_i\,,a_j] (a_1 \cdots a_n)(b_1 \cdots b_n)\:,\\
\label{nil2formel2} d_{\mu_{n-1}}(a,b)   &=& \prod_{1\le i < j \le n} [b_i\,,a_j]\:.
\end{eqnarray}
In fact, if $i<j$, we have $b_ia_j = [b_i,a_j]a_jb_i$\,; in shuffling all the factors $b_i$ in the product $a_1b_1 \cdots a_nb_n$
to the right, one introduces all the commutators $[b_i\,,a_j]$ with $i<j$. But these are central in $G$ so can be gathered on the
left which proves the first formula. It implies the second one as \hfill $d_{\mu_{n-1}}(a,b) = (a_1b_1 \cdots a_nb_n)(b_1\cdots
b_n)^{-1} $\\$ (a_1\cdots a_n)^{-1}$.

Now we prove the desired equivalences.\V

\N(1) $\Rightarrow$ (2). We have $D_{\mu_{n-1}} \subset G\st \subset Z(G)$ by equation \REF{nil2formel2} and as $G$ is 2-step
nilpotent. Each commutator $[b_i,a_j]$ being bilinear identity \REF{nil2formel2} shows that $d_{\mu_{n-1}}$ is bilinear, too.\V

\N(2) $\Rightarrow$ (3). The map $K\to G^n$, $x\mapsto (f_1(x)\,,\ldots, f_n(x))$, is linear. Hence by Proposition \ref{quadolin} below the
composite map with the quadratic map \map{\mu_{n-1}}{G^n}{G} is still quadratic.\V

\N(3) $\Rightarrow$ (4). It suffices to take two of the maps $f_i$ equal to the identity of $G$ and the others equal to the trivial
map.\V

\N(4) $\Rightarrow$ (1). Using  \REF{df(a,b)=}  we get $d_{2_G}(a,b)=[-a,b]$, hence $[-a,b] \cdot{}^b[-a,b\st] = [-a,bb\st] = 
d_{2_G}(a,bb\st) = d_{2_G}(a,b) d_{2_G}(a,b\st) = [-a,b] [-a,b\st]$, for $a,b,b\st\in G$.  It follows that ${}^b[-a,b\st] = 
[-a,b\st]$, whence $G\st$ is central in $G$ and $G$ is 2-step nilpotent.\V

\N(5) It suffices to take $K=G$ and  each $f_i=id$ in assertion (3).\hfill $\Box$\V




We now recall some relations already proved in Lemma 2 in \cite{Niq}.

\begin{prop}\label{Propquad}   Let \map{f}{G}{H} be a \qu\ map.  Then the following relations
hold for $a,b\in G$.

\comment{
\BE\label{f(-a)}    f(a^{-1})  = {} -f(a) + d_f(a,a)  \EE
\BE\label{f(a-b)}   f(ab^{-1}) = f(a) - f(b) - d_f(ab^{-1},b)  \EE
\BE\label{f[a,b]}   f[a,b] = [f(a),f(b)] + d_f(a,b) - d_f(b,a)  \EE
\BE\label{d-f}      d_{-f}(a,b) = {} - d_f(b,a) - f[a,b]  \EE
}

\comment{
\begin{eqnarray}
 \label{f(-a)}     f(a^{-1})   &=&  {} -f(a) + d_f(a,a)  \\
 \label{f(a-b)}   f(ab^{-1})  &=&  f(a) - f(b) - d_f(ab^{-1},b) \\
 \label{f[a,b]}   f[a,b]  &=&  [f(a),f(b)] + d_f(a,b) - d_f(b,a)  \\
 \label{d-f}      d_{-f}(a,b)  &=&  {} - d_f(b,a) - f[a,b]   
\end{eqnarray}
}

\begin{eqnarray}
 \label{f(-a)}    \makebox[15mm][l]{$f(a^{-1})$}   &=&  {} -f(a) + d_f(a,a)  \ruled \\
 \label{f(a-b)}   \makebox[15mm][l]{$f(ab^{-1})$}  &=&  f(a) - f(b) - d_f(ab^{-1},b) \ruled \\
 \label{f[a,b]}   \makebox[15mm][l]{$f[a,b]$}  &=&  [f(a),f(b)] + d_f(a,b) - d_f(b,a)  \ruled \\
 \label{f(a+b-a)}\makebox[15mm][l]{$f({}^ab)$}  &=&  {}^{f(a)} f(b)  + d_f(a,b) - d_f(b,a)  \ruled \\
 \label{d-f}      \makebox[15mm][l]{$d_{-f}(a,b)$}  &=&  {} - d_f(b,a) - f[a,b]   
\end{eqnarray}

\end{prop}

We point out that relation \REF{f[a,b]} has a conceptual meaning which makes it a crucial ingredient in the structure theory in section 2; it also is the only relation which generalizes to the setting of semi-abelian categories where it again plays a key role  \cite{QuadFunct}.\V


\N{\bf \large Composition of \qu\ maps.}\V

\N The content of this section is essentially due to my former student O.\ Perriquet \cite{Pe}.
Composing two \qu\ maps does not give a \qu\ map in general, see Example \ref{Exqumaps}(0).
However, the following  sufficient
condition assuring stability under composition is often satisfied in practice, see \cite{BHP}, \cite{ModQuad}.

\begin{defi}\label{qupairdef} \rm Let $K \maprtop{g} G \maprtop{f} H$ be two \qu\ maps relative some subgroup $A$ of $K$ and $B$ of $G$, resp. We
say that the couple $(f,g)$ is a {\em \qu\ pair}\/ if $g(A) \subset B$ and $d_f(D_g \times G) = d_f(G \times D_g) = 0$.
\end{defi}

For example, if $g$ and $f$ are \qu\ maps such that $D_g\subset G'$ then for $A=K'$ and $B=G'$, $(f,g)$ is a \qu\ pair by \REF{f[a,b]}.

\begin{prop}\label{qupairprops} If $(f,g)$ as above is a  \qu\ pair  the composite map $f\rond g$ is \qu\ relative $A$, with 
  \BE\label{dfrondg}   d_{f\rond g} = f_{\raisebox{-1mm}{$*$}} d_g + (g \times g)^* d_f\:.   \EE
Moreover, if $D_g \subset BG\st$, $g$ induces a linear map
\map{\bar{g}}{K/AK\st}{G/BG\st}  and we have
  \BE\label{wfrondg}  w_{f\rond g} = f_{\raisebox{-1mm}{$*$}} w_g + (\bar{g} \otimes \bar{g})^* w_f\:.  \EE
\end{prop}
 
\N We call \REF{dfrondg} or \REF{wfrondg} the {\em derivation property}\/ of a \qu\ pair.\V

\N{\bf Proof\,:} Writing $G$ and $H$ additively we have for $a,b\in K$
\begin{eqnarray*}
d_{f\rond g}(a,b)  &=& fg(ab) -fg(b) -fg(a) \\
&=&  f(d_g(a,b) + g(a) + g(b)) - fg(b) -fg(a) \\
&=&  f(d_g(a,b)) + f(g(a) + g(b)) - fg(b) -fg(a) \quad \mbox{since $d_f(D_g \times G) =0$}\\
&=&  f(d_g(a,b)) + d_f(g(a) , g(b)) + fg(a) +  fg(b) - fg(b) -fg(a) \\
&=&  f(d_g(a,b)) + d_f(g(a) , g(b)) \:.
\end{eqnarray*}
Next we check that $d_{f\rond g}$ is linear in the first variable; the argument for the second variable is similar. For $a\st \in
K$, we have 
\begin{eqnarray*}
d_{f\rond g}(aa\st,b)  &=& f(d_g(aa\st,b)) + d_f(g(aa\st) , g(b)) \\
 &=&  f(d_g(a ,b) +  d_g( a\st,b)) +  d_f(g(a ) + g( a\st) + d_g(a,a\st) , g(b)) \\
&\stackrel{(*)}{=}&  f(d_g(a ,b)) +  f(d_g( a\st,b)) +  d_f(g(a ),g(b))  + d_f(g( a\st) , g(b)) \\ 
&\stackrel{(**)}{=}&  f(d_g(a ,b)) +  d_f(g(a ),g(b))  + f(d_g( a\st,b)) +   d_f(g( a\st) , g(b)) \\
&=&   d_{f\rond g}(a,b) + 
d_{f\rond g}(a\st,b)\:.
\end{eqnarray*} 
Equations $(*)$ and $(**)$ hold since $d_f(D_g \times G) =0$ and since $D_f\subset
Z(I_f)$, resp.

Moreover, if $a$ or $b$ is in $A$ then $d_{f\rond g}(a,b)=0$ since then $d_g(a,b) =0$ and $d_f(g(a) , g(b))=0$ as $g(A)\subset B$.
It remains to check that $D_{f\rond g} \subset Z(I_{f\rond g})$. Let $a,b,c \in K$. Then $[d_f(g(a ),g(b)) , f(g(c))] = 0$ since
$D_f \subset Z(I_f)$, and 
  \[ [f(d_g(a ,b)) , f(g(c))] =  f[d_g(a ,b),g(c)] - d_f(d_g(a ,b),g(c)) + d_f(g(c) , d_g(a ,b)) \]
by \REF{f[a,b]}. But the first of the latter three terms is trivial since $D_g \subset Z(I_g)$, the other two since $d_f(D_g \times
G) = d_f(G \times D_g) =0$. \hfill $\Box$\V

\begin{kor}\label{quadolin} Pre- or postcomposing a quadratic map by a linear map gives a \qu\ map. More precisely, if $f\,\colon
G\to H$ is
\qu\ and $g\,\colon K\to G$, $h\,\colon H\to L$ are linear maps of groups  then $hfg$ is \qu\ with $d_{hfg} = h\lstar (g\times g)^*
d_f$.\hfill$\Box$
\end{kor}\vspace{2mm}

\begin{kor}\label{Quadfunct} There is a functor  ${\rm Quad}(G,-)$ from the category of groups to the category of sets sending a group $H$ to ${\rm Quad}(G,H)$ and a homomorphism $f\mapco H \to K$ to the map $f_*\mapco {\rm Quad}(G,H)\to {\rm Quad}(G,K)$.
\end{kor}\vspace{2mm}

The following relations are crucial in dealing with quadratic categories, see \cite{BHP}.

\begin{prop}\label{2Hfqupair} Let \map{f,g}{G}{H} be two functions between groups.\V

\N(1) If $(2_H,f)$ is a quadratic pair then $D_f$ is central in $H$.\V

\N(2) If $(2_H,g)$ is a quadratic pair then the following relations holds for $a,b \in G$.
\begin{eqnarray*} d_{f+g}(a,b)  &=& d_f(a,b) +d_g(a,b) + d_{2_H}(g(a)\,,f(b)) \\
  &=& d_f(a,b) +d_g(a,b) + [f(b)\,,g(a)]\:.
\end{eqnarray*}
\end{prop}

\N{\bf Proof\,:} If $2_H$ is \qu\ then $d_{2_H}(a,b) = [-a,b]$ by \REF{df(a,b)=}, whence (1) follows from the equations $[H,D_f]
= [-H,D_f] = d_{2_H}(H\times D_f) = 0$. To prove (2) calculate
\begin{eqnarray*}
d_{f+g}(a,b)  &=& (f+g)(a+b) - (f+g)(b) - (f+g)(a) \\
  &=&  f(a+b) + g(a+b) - g(b) - f(b) - g(a) - f(a) \\
  &=& d_f(a,b) + f(a) + f(b) + d_g(a,b) + g(a) + g(b) - g(b) - f(b) - g(a) - f(a) \\
  &\stackrel{(*)}{=}& d_f(a,b)+ d_g(a,b)  + f(a) + f(b) + g(a) - g(a)  - f(b) + [f(b)\,,g(a)] - f(a) \\
  &=& d_f(a,b) +d_g(a,b) + [f(b)\,,g(a)]\:.
\end{eqnarray*}
The last equation follows from the fact that $H\st$ is central in $H$ as $H$ is 2-step nilpotent since $2_H$ is quadratic, see Proposition \ref{quadnil2}.
Equation $(*)$ is due to the fact that $D_g$ is central in $H$ by assertion (1). Finally, as $H$ is 2-step nilpotent,
$[f(b)\,,g(a)]= -[g(a)\,,f(b)] = [-g(a)\,,f(b)] = d_{2_H}(g(a)\,,f(b))$.\hfill$\Box$\vspace{8mm}

\N{\bf \large The category of quadratic pairs.}\V


\begin{defi}\label{qupairmapdef}  A pair of groups $(G,B)$ consists of a group $G$ together with a subgroup $B$ of $G$.  A  pair map $f\mapco 
(G,B) \to (H,C)$ between pairs of groups  is a function $f$ from $G$ to $H$ such that $f(B)\subset C$. Moreover, a pair map $f$ is a linear pair map if the function $f$ is linear, and $f$ is a quadratic pair map if $f$ is quadratic relative $B$ such that $C$ contains $D_f$.
\end{defi}\V

Note that any quadratic map $f\mapco G \to H$ is a quadratic pair  map from $(G,G')$ to $(H,H'D_f)$ by \REF{f[a,b]} and from $(G,\rad{f})$ to $(H,f(\rad{f})D_f)$ by Proposition \ref{qumapsprops}(c).

\begin{prop}\label{qupairmapcomp} Let $(K,A) \mr{g} (G,B) \mr{f} (H,C)$ be quadratic pair maps. Then the composite map $f g \mapco (K,A) \to (H,C)$
is a quadratic pair map whose deviation satisfies the derivation rules \REF{dfrondg} and \REF{wfrondg}.
\end{prop}

This is immediate from Proposition \ref{qupairprops}; it leads to the following generalization of the ``quadratic envelope  of the category of 2-step nilpotent groups"  constructed by Jibladze and Pirashvili \cite{Niq}: this is the category, denoted by {\bf Niq}, whose objects are 2-step nilpotent groups and whose morphisms from $G$ to $H$ are the \qu\ maps $f\mapco G \to H$ such that $D_f \subset H'$.\V

\begin{kor}\label{cat QP} Pairs of groups and quadratic pair maps between them form a category denoted by ${\bf QP}$ which we call 
quadratic envelope  of the usual category of linear pair maps.
\end{kor}

In fact, the category {\bf Niq} fully embeds into ${\bf QP}$ by sending $G$ to the pair $(G,G')$.

Denote by ${\bf Gr}$ and ${\bf Ab}$ the category of groups and abelian groups, resp. Then the following is again an immediate consequence of Proposition \ref{qupairprops}.\V

\begin{prop}\label{qupairmapcats} Let ${\bf NQP}$ resp.\ ${\bf AQP}$ resp.\ ${\bf CP}$ be the full subcategory of ${\bf QP}$ consisting of those pairs of groups $(G,B)$ for which $B$ is normal in $G$, resp.\  abelian, resp.\ central containing $G'$.\V

 (a) There are functors 
  $ Gr  \mapco {\bf  QP} \to {\bf Ab} \times {\bf Gr} \:,\quad Gr_N \mapco {\bf NQP} \to {\bf Gr} \times {\bf Gr} \:,\quad $ $
  Gr_A \mapco {\bf AQP} \to {\bf Ab} \times {\bf Ab} \:,\quad Gr_C \mapco {\bf CP} \to {\bf Ab} \times {\bf Ab} $ such that 
 $Gr$ and $Gr_A$ send  an object $(G,B)$   to the object $(G/G'B,B)$  while $Gr_N$ and $Gr_C$ send $(G,B)$ to $(G/B,B)$, and all of them send $f\mapco 
(G,B) \to (H,C)$ to $(\bar{f},f_B)$ where $f_B\mapco B \to C$ is the restriction of $f$ and $\bar{f}$ is induced by $f$.\V


  (b) There is  a bifunctor  $D\mapco {\bf AQP}^{op} \times {\bf AQP} \hspace{1mm}\to\hspace{1mm} {\bf Ab}$ defined by \[ D((G,B),(H,C)) = {\rm Hom}((G/G'B) \ot (G/G'B) \,, C)\quad\mbox{and}\quad D(g^{op},f) = f_*(\bar{g}\ot \bar{g})^*\:.\]
  
  (c) Assigning to a map $f\mapco 
(G,B) \to (H,C)$ in ${\bf AQP}$ its defect $w_f \in D((G,B),$ $(H,C))$ defines a {\em derivation}\/ from ${\bf AQP}$ to $D$, see \cite{BW}.

\end{prop}

Note that  ${\bf CP} $ is contained in the intersection of ${\bf NQP} $ and ${\bf AQP} $, and that for $(G,B)\in  {\bf CP} $ the group $G$ is 2-step nilpotent. Thus  the category {\bf Niq} also identifies with the full subcategory of {\bf CP} consisting of the objects $(G,G')$.
Moreover, ${\bf CP} $ is a right quadratic category in the sense of \cite{BHP}; from several points of view, it plays the same role in quadratic algebra as  the category of abelian groups plays in classical algebra, see \cite{ModQuad}; it can therefore be considered as a different generalization of the category of abelian groups than the category of square groups constructed in \cite{BP}. The relation between these two generalizations, however, is not yet understood.
\vspace{6mm}

\N{\bf \large Universal relative quadratic map.}\V

We now show that the functor ${\rm Quad}(G,-)$ is representable; this result is also obtained in \cite{Niq}. We thus get an endofunctor $Q$ of the category of groups the properties of which are studied.\V

\begin{satz}\label{Q(G)} Let $G$ be a group and $B$ a subgroup of $G$.\V

  (i) There exists a universal \qu\ map relative $B$, \map{q=q_{G,B}}{G}{Q(G,B)} to some group $Q(G,B)$, i.e., for any \qu\ map of
groups 
$f\,\colon G\to H$ relative
$B$ there exists a unique linear map $\hat{f}\,\colon Q(G,B) \to H$ such that $\hat{f} q = f$.\V

  (ii) The sequence of group homomorphisms 
\BE\label{Qsequ}  0 \:\lra\: \GBoGB \hspace{2mm} \maprtop{w_q} \hspace{2mm}  Q(G,B) \hspace{2mm} \maprtop{\widehat{id}} \hspace{2mm} G \:\lra\:
1 \EE
 is exact. Actually, it is a central group extension represented by the bilinear 2-cocycle $D\,\colon G\times G \to \GBoGB$,
$D(a,b) = {} - \bar{a} \ot \bar{b}$ which corresponds to the canonical section $q$ of $\widehat{id}$.
\end{satz}

\N{\bf Proof\,:} The trick is to {\em define}\/ the group $Q(G,B)$ by the cocycle $D$, i.e., we let $Q(G,B) = (\GBoGB) \times G$
endowed with the group law $(x,a) + (y,b) = (x+y - \bar{a} \ot \bar{b}\,, ab)$. Furthermore, let $q\,\colon G\to Q(G,B)$, $q(a) =
(0,a)$. Then for $a,b\in G$, $q(a)+q(b) = ( - \bar{a} \ot \bar{b}\,,ab)  =  ( - \bar{a} \ot \bar{b}\,, 1) + q(ab)$, whence
$d_q(a,b) = (\bar{a} \ot \bar{b}\,, 1)$. This term is bilinear, central in $Q(G,B)$, and trivial whenever $a$ or $b$ is in $B$, so
$q$ is a \qu\ map relative $B$ with $w_q(x) = (x,1)$. In order to prove its universal property let $f\,\colon G\to H$ be some \qu\
map relative
$B$. Define
$\hat{f}\,\colon Q(G,B) \to H$ by $\hat{f} (x,a) = w_f(x) + f(a)$. Then $\hat{f}$ satisfies $\hat{f} q = f$, and for $(x,a),(y,b)
\in Q(G,B)$, we have 
  \begin{eqnarray*}
\hat{f} ((x,a) + (y,b)) &=& w_f(x+y-\bar{a} \ot \bar{b}) + f(ab) \\
 &=& w_f(x) + w_f(y) - d_f(a,b) + d_f(a,b) + f(a) + f(b) \\
 &=& w_f(x) + f(a) + w_f(y)  + f(b) \hspace{5mm}\mbox{since $\Imm{w_f} = D_f \subset Z(I_f)$}\\
  &=& \hat{f}(x,a)  +  \hat{f}(y,b)\:.
   \end{eqnarray*}
To prove uniqueness of $\hat{f}$ let \map{g}{Q(G,B)}{H}  be a linear map such that $gq =f$. Then $g(0,a) = gq(a)=f(a)$ and 
$g(\bar{a} \ot \bar{b},1) = g\,d_q(a,b) = d_{gq}(a,b)$ (by \ref{quadolin}) $= d_f(a,b) = w_f(\bar{a} \ot \bar{b})$, whence by
linearity, $g(x,1) = w_f(x)$ for all $x\in \GBoGB$. Thus $g(x,a) = g((x,1)+(0,a)) = w_f(x) + f(a) = \hat{f}(x,a)$, whence
$g=\hat{f}$. Finally, we have $\widehat{id}(x,a) = w_{id}(x)a = a$ as $id$ is linear, so $\widehat{id}$ is the projection to the
second factor which proves exactness of the sequence in (ii) as we saw already that $w_q$ is the canonical injection of the first
factor.\hfill $\Box$\vspace{3mm}

\nc{\res}[2]{{#1}_{| \raisebox{-0.7mm}{$\sst #2$}}}

Write $Q(G) = Q(G,\{1\}) = Q(G,G')$ and $q=q_G=q_{G,\{1\}}$. Note that Theorem \ref{Q(G)}(i) says that the map 
  \BE\label{Qrepres} q^*\mapco {\rm Hom}(Q(G),H) \to {\rm Quad}(G,H)\EE
   is a bijection natural in $H$.\V

\begin{kor}\label{Q(free)} Let $X$ be a set and $F$ a free group with basis $X$. Then there is a (non natural) isomorphism $\phi \mapco Q(F) \to (F\ab \otimes F\ab) \times F$ such that $\phi q_F$ is given by
\BE\label{phiqF} \phi q_F ( \prod_{i=1}^n x_i^{\epsilon_i}) = \Big( \sum_{i=1}^n \frac{1-\epsilon_i}{2} \,\bar{x}_i \ot \bar{x}_i + \sum_{1\le i<j\le n} \epsilon_i \epsilon_j \,
\bar{x}_i \ot \bar{x}_j \,, \prod_{i=1}^n x_i^{\epsilon_i} \Big)\EE
for $x_i\in X$, $\epsilon_i=\pm 1$. Consequently, for any group $H$ there is a bijection
\BE\label{Quad(F,H)} Quad(F,H) \hcong \{(\chi,\psi) \in H^X \times H^{X\times X} \,|\, [\Imm{\chi},\Imm{\chi}] =  [\Imm{\chi},\Imm{\psi}]=\{1\}        \}\EE
which carries $f\in Quad(F,H)$ to $(\res{f}{X}\,,\res{d_f}{X\times X})$.
\end{kor}

\proof The map $\res{q_F}{X}$ induces a splitting of the central extension \REF{Qsequ}. The formula for $\phi q_F$ is obtained using \REF{f(ab)=} and \REF{f(-a)}. To determine ${\rm Quad}(F,H)$ now use the bijection $q_F^*$ in \REF{Qrepres} and the fact that $X\times X \to F\ab \ot F\ab$, $(x,y)\mapsto x\ot y$, is a basis of the abelian group $F\ab \ot F\ab$.\hfbox\V

\begin{prop}\label{Qfunctor} There is an endofunctor $Q$ of the category of groups sending $G$ to $Q(G)$ and $f\mapco G \to H$ to $Q(f) = \widehat{q_Hf}\mapco Q(G) \to Q(H)$. It has the following properties.\V

(a) If $G$ is $n$-step  nilpotent for $n\ge 2$ then so is $Q(G)$.\V

(b) If $G_1 \mr{\alpha} G_2 \mr{\beta} G_3 \to 1$ is an exact sequence of groups then so is the sequence 
\BE\label{Qexsequ} Q(G_1) \hspace{1mm}\times\hspace{1mm} G_1\ab \ot G_2\ab \hspace{1mm}\times\hspace{1mm} G_2\ab \ot G_1\ab \hmr{\xi} Q(G_2) \hmr{Q(\beta)}  Q(G_3) \hspace{1mm}\to \hspace{1mm} 1\EE
where $\xi(x,y,z) = Q(\alpha)(x) + w_q(\alpha\ab \ot 1)(y) + w_q(1\ot \alpha\ab)(z)$. One also has the identity 
\BE\label{KerQ(f)}  \Ker{Q(\beta)} = \Imm{q_{G_2}\alpha} + \Imm{w_{q_{G_2}}(\alpha\ot \alpha + \alpha\ab\ot 1_{G_2} + 1_{G_2} \ot \alpha\ab)} \EE

\end{prop}

Similarly, assigning the group $Q(G,B)$ to a pair $(G,B)$ defines a functor from the category of linear pair maps to the category of groups.\V

\proof (a): Let $a\in Q(G)$, $b\in \gamma_n(Q(G))$; we must show that $[a,b]=1$. Modulo the central subgroup \Imm{w_q} we may assume that $a=qa',b=qb'$ for some $a',b'\in G$. Then by \REF{f[a,b]}, $[a,b] = q[a',b'] - w_q(\overline{a'} \ot \overline{b'} - \overline{b'} \ot \overline{a'}) =1$ since 
$[a',b']=1$ as $b' \in \gamma_n(G)$ and since $\overline{b'}=0$ as $n\ge 2$.

For (b) use naturality of   sequence \REF{Qsequ} with respect to the maps $\alpha$ and $\beta$; then exactness of sequence \REF{Qexsequ} follows by an easy diagram chase together with right exactness of the tensor product. Identity \REF{KerQ(f)} then follows using the fact that $Q(G_1) = \langle \Imm{q} \rangle + \Imm{w_q}$.\hfbox

\begin{kor}\label{Qgen} If $G$ is generated by a subset $X$ then $Q(G)$ is generated by the subset $q(X) \cup d_q(X\times X)$.
\end{kor}
Just apply Proposition \ref{Qfunctor}(b) to the epimorphism from the free group of basis $X$ to $G$ and use Corollary \ref{Q(free)}.\V

\begin{kor}\label{f=g} Two quadratic maps $f,g\mapco G\to H$ coincide if and only if they coincide on a generating subset $X$ of $G$ and $d_f$ coincides with $d_g$ on $X\times X$.\end{kor}

Just note that $f=g$ iff $\hat{f}=\hat{g}$ as $q$ is injective (having $\widehat{id}$ as a retraction).\V

The above properties of the functor $Q$ allow to construct quadratic maps on groups defined by generators and relations, as follows. 
Let $G=\langle X|R \rangle$ be a presentation of $G$, i.e.\ let $X$ be a set, $F$ the free group with basis $X$, $R$ a subset of $F$ and  $\pi\mapco F\auf G$ a homomorphism whose kernel is the normal subgroup of $F$ generated by $R$. For $r\in R$, write $r =\prod_{i=1}^{n_r} x_{ri}^{\epsilon_{ri}}$ with $ x_{ri}\in X$ and $\epsilon_{ri} =\pm 1$, and $\bar{r} =rF'= \sum_{x\in X} k_{rx}\bar{x}$ in $F\ab$, with $k_{rx}\in \Z$.\V

\begin{prop}\label{qumapsG=XR} For some group $H$ let $(\chi,\psi) \in H^X \times H^{X\times X}$. Then there exists a quadratic map $f\mapco G\to H$ such that 
\BE\label{chipsi} \mbox{$f\pi(x) = \chi(x)$ \hspace{3mm} and\hspace{3mm}  $d_f(\pi(x),\pi(y)) = \psi(x,y)$ \hspace{3mm} for\hspace{1mm}  $x,y\in X$} \EE
if and only if for all $r\in R$ and $y\in X$ the following three conditions hold.
\begin{rom}

\item \hfill $[\Imm{\chi},\Imm{\chi}] =  [\Imm{\chi},\Imm{\psi}]=\{1\}$ \hfill\phantom{i}

\item \hfill $\sum_{i=1}^{n_r} \Big(  {\epsilon_{ri}}\chi(x_{ri}) +  \frac{1-\epsilon_i}{2} \,\psi(x_{ri}\,,x_{ri}) \Big) + \sum_{1\le i<j\le n} \epsilon_{ri} \epsilon_{rj} \,
\psi(x_{ri}\,,x_{rj}) \hspace{1mm}=\hspace{1mm} 0$  \hfill\phantom{i}

\item \hfill $\sum_{x\in X} k_{rx}\psi(x,y) \hspace{1mm}=\hspace{1mm}  \sum_{x\in X} k_{rx}\psi(y,x)
\hspace{1mm}=\hspace{1mm} 0$  \hfill\phantom{i}
\end{rom}

\N Moreover, any \qu\ map from $G$ to $H$ is induced by maps $\chi$ and $\psi$ as above.

\end{prop}\V

\proof Let $F_1$ be the free group with basis $F\times R$. Then  the sequence $F_1\mr{\delta} F\mr{\pi} G \to 1$ is exact  with $\delta(a,r)={}^ar$ for $(a,r)\in F\times R$.
Hence by \REF{Qrepres} and Proposition \ref{Qfunctor}(b) there exists a \qu\ map $f\mapco G\to H$ satisfying \REF{chipsi} iff there is a homomorphism $\kappa\mapco Q(F) \to H$ factoring through $Q(\pi)$ and satisfying the property
 \BE\label{kappa} \mbox{$\kappa q_F(x) = \chi(x)$ \hspace{3mm} and\hspace{3mm}  $\kappa d_{q_F}(x,y)  = \psi(x,y)$ \hspace{3mm} for\hspace{1mm}  $x,y\in X$\,.} \EE
 By Corollary \ref{Q(free)} a homomorphism $\kappa$ satisfying \REF{kappa} exists iff  condition (i) holds, so suppose this true in the sequel. By Proposition \ref{Qfunctor}(b) $\kappa$ factors through $Q(\pi)$ iff $\kappa \,Q(\delta) = \kappa w_{q_F} (\delta\ab \ot 1_{F\ab} + 1_{F\ab} \ot \delta) = 0$. But $\Imm{\delta\ab} = \langle \Imm{R\to F\ab} \rangle$ whence the second identity holds iff for all $(r,x)\in R\times X$, $\kappa w_{q_F}(\bar{r},\bar{x}) = \kappa w_{q_F}(\bar{x},\bar{r}) =0$ which by expanding $\bar{r}$ is equivalent to condition (iii). Next consider $\kappa \,Q(\delta)$. By Corollary \ref{Qgen} $Q(F_1)$ is generated by $q_F\delta(F\times R) + \Imm{w_{q_F}(\delta\ab \ot\delta\ab)}$. Let $(a,r)\in F\times R$. Then $q_F\delta(a,r) = q_F({}^a r) = {}^{q_F(a)}q_F(r) + 
 w_{q_F}(\bar{a}\ot \bar{r} - \bar{r}\ot \bar{a})$ by \REF{f(a+b-a)}. Hence if $\kappa$ satisfies conditions (iii) it annihilates $\Imm{Q(\delta)}$ iff it annihilates $q_F(R)$ which is equivalent to condition (ii) by \REF{phiqF}.\hfbox\V

\vspace{3mm}

\section{Relation with Passi's construction.}

In this section we study relative polynomial maps in the sense of Passi by using the nonabelian theory of the first section as an
essential tool. We will see that the proof of the main properties in the quadratic case becomes more natural in this way than
in  former approaches in the literature.

For a commutative ring $R$ with unit  let $I_R(G)$ denote the augmentation ideal of the
group algebra $R(G)$; for $k\ge 0$, $I_R^k(G)$ denotes its $k$-th power, with the convention  $I_R^0(G) = R(G)$. If $R=\Z$ we also
write $I(G)=I_{\Z}(G)$. For a function $f\,\colon G\to A$ to some {\em
abelian}\/ group $A$ let $\bar{f}\,\colon \Z(G) \to A$  denote the  extension of $f$ to a \Z-linear
homomorphism.

\begin{defi} Let $G$ be a group and $B$ be a normal subgroup of $G$. We say that a function $f\,\colon G\to A$ as above is   polynomial of degree $\le n$ relative $B$ if 
$\bar{f}$ annihilates the subset $1+I(B)I(G)+I^{n+1}(G)$ of $\Z(G)$. 

Moreover, we say that $f$ is (normalized) polynomial of degree $\le n$ if it is polynomial of degree $\le n$
relative $\{1\}$. 
\end{defi}

\begin{rem} \rm The notion of (absolute) polynomial map from groups to {\em abelian}\/ groups is due to Passi \cite{Pa68}. The
relative case was only implicitely present in most of the work in the literature based on this notion, notably when  related to  the
dimension subgroup problem. See \cite{Pa} for a thorough treatment of the subject.
\end{rem}

 Note that $f$ is polynomial of degree $\le 0$ iff $f=0$. Moreover, if $f$ is polynomial of degree $\le n$, $f$ is also polynomial of degree $\le n$ relative
$\gamma_n(G)$ as
$I(\gamma_n(G)) \subset I^n(G)$; this is immediate by inductive application of the formula
 \BE\label{comm-1} [a,b] - 1 = [a-1\,,b-1]a^{-1}b^{-1} \EE
for $a,b\in G$ where $[-,-]$ denotes the group commutator on the left and the ring commutator on the right.\V

The following inductive characterization of polynomial maps is useful for proving polynomiality, see \cite{Q3} where a more general
version is developed (with respect to arbitrary $N$-series). For convenience of the reader we give a direct proof of our special
case here which is very short and easy anyway.\vspace{2mm}

\begin{prop}\label{polind} Let \map{f}{G}{A} be any normalized function from a group $G$ to some abelian group $A$. Then we have the following properties.\V

\N(1) For $a,b \in G$, $d_f(a,b) = \bar{f}((a-1)(b-1))$.\V

\N(2) For $n\ge 1$, $f$ is polynomial of degree $\le n$ relative some given normal subgroup $B$ of $G$ if and only if the following two conditions
hold.\V

  (a) The map \map{d_f(a,-)}{G}{A} (or equivalently, \map{d_f(-,a)}{G}{A}) is polynomial of degree $\le n-1$ for all $a\in G$.\V

  (b) For all $(b,a)\in B\times G$, $f(ba) = f(b) + f(a)$ or equivalently, $d_f(b,a)=0$.
\end{prop}\vspace{2mm}

\N{\bf Proof\,:} (1) is immediate from expanding the right hand term. To prove (2) let $a\in G$. Then for $b\in G$,
  $\overline{d_f(a,-)}(b-1) = d_f(a,b) - d_f(a,1) = \bar{f}((a-1)(b-1)) $ by (1). As the elements $b-1$, $b\in G$, generate $I(G)$
as a \Z-module, it follows by linearity that $\overline{d_f(a,-)}(x) = \bar{f}((a-1)x) $ for all $x\in I(G)$. Hence 
$\overline{d_f(a,-)}(I^n(G)) = \bar{f}((a-1) I^n(G)) $, which implies that property (a) is equivalent to $f$ being polynomial of
degree $\le n$. Moreover, (b) is equivalent to $\bar{f}((b-1)(a-1)) = d_f(b,a) = 0$ for all $(b,a)\in B\times G$  which in turn 
means that $\bar{f}(I(B)I(G)) = 0$. \hfill $\Box$\vspace{4mm}

 In low degrees, we obtain the following caracterization of (relative) polynomial maps.

\begin{kor}\label{pollowdeg} Let  $f\,\colon G\to A$ and $B$ as in \ref{polind}.\V


\N(1) $f$ is polynomial of degree $1$ relative $B$ iff it is linear.\V

\N(2) $f$ is polynomial of degree $2$ relative $B$ iff $d_f$ is bilinear and annihilates $B$ in the first variable.\hfill $\Box$
\end{kor}

\begin{kor}\label{quad=pol} Let $G$ be a group, $B$ be a {\em central} subgroup of $G$ and $f\,\colon G\to A$ as in \ref{polind}.
Then $f$ is polynomial of degree $\le 2$ relative $B$ if and only if $f$ is \qu\ relative $B$.\end{kor}

This is immediate from \ref{pollowdeg}, just note that centrality of $B$ implies that $d_f(a,b) = \bar{f}((a-1)(b-1))
=\bar{f}((b-1)(a-1)) = d_f(b,a)$ for $(a,b)\in G\times B$. \hfill $\Box$

Before exploiting corollary \ref{quad=pol} we give some examples the verification of which is based on \ref{polind} and
\ref{pollowdeg}.

\begin{bspe}\rm Let $G$ be a group.\V

\N(1) If $G$ is 3-step nilpotent, the commutator map $G\times G \to G$, $(a,b)\mapsto [a,b]$, is bipolynomial of degree $\le 2$.\V

\N(2) Recall that the non-abelian tensor square $ G\ot G$ of $G$ is a group closely related to the homotopy group $\pi_3\Sigma K(G,1)$ and also to the second homology group $H_2(G)$, see \cite{BrLo} and \cite{Ellis}. 
Now if $G$ is 2-step nilpotent, the natural map $G\times G \to G\ot G$, $(a,b)\mapsto a\ot b$, is bipolynomial of degree $\le
2$, see \cite{GoG} where this fact is used to compute $G\ot G$ for 2-step nilpotent groups.\V

\N(3) For $G$ abelian the $n$-fold diagonal map \map{\delta^n}{G}{G\htt{n}}, $\delta^n(a) =  {a} \ot \cdots \ot  {a}$, is 
polynomial of degree $\le n$, see \cite{Q3}.

\end{bspe}

In order to introduce {\em universal}\/ relative polynomial maps we recall the following definition and facts from \cite{PolProp}.\V

\begin{defi}\label{defi} \rm \quad Let $G$ be a group, $R$ as above, and
$B \lhd G$  a normal subgroup. We define the quotient $R$-algebra without unit
   \[ P_{n,R}(G,B) = I_R(G)\Big/ \left(I_R(B) I_R(G) +
I_R^{n+1}(G) \right) \:. \]
The canonical quotient map
from
$I_R(G)$ to
$P_{n,R}(G)$ or to $P_{n,R}(G,B)$ is denoted by $\rho$ or $\rho_n$. If $f\mapco (G,B) \to (H,C)$ is a linear pair map in ${\bf NQP}$ it induces a morphism of $R$-algebras $P_{n,R}(f) \mapco P_{n,R}(G,B) \to P_{n,R}(H,C)$ defined by $P_{n,R}(f)\rho(a-1) = \rho(f(a)-1)$ for $a\in G$.

We point out that  $\:P_{n,R}(G) = P_{n,R}(G,\{1\})\:$ is the polynomial group constructed by
Passi \cite{Pa68}; the relative version was introduced and studied modulo torsion in \cite{PolProp}. The aim of this section is to
explore in detail the structure of  $P_{2,\Z}(G,B)$ for central $B$ which is the crucial ingredient of our abelian models for 2-step
nilpotent groups in \cite{Diss}, \cite{Habil}.

\end{defi}

Using the elementary identification 
   \BE\label{Iquot} I_R(G)/I_R(B)R(G) \:\stackrel{\cong}{\to}\: I_R(G/B)\:,\quad \mbox{$\overline{a-1} \mapsto
\bar{a}-1\hspace{3mm}$ for $ a\in G$,} \EE
 we see  that multiplication in the ring $P_{n,R}(G,B)$ gives rise to an $R$-linear map
\BE\label{mundef} \mu_n\:\colon\: P_{n-1,R}(G/B) \otimes_{\raisebox{-1mm}{$\sst R(G)$}} P_{n-1,R}(G,B) \:\lra \: P_{n,R}(G,B)  \EE
such that for $x,y\in I_R(G)$, $\mu_n(\rho_{n-1}(x) \otimes \rho_{n-1}(y)) = \rho_n(xy)$.
This shows that via left multiplication, $P_{n,R}(G,B)$ is a  nilpotent left $R(G/B)$-module of class $\le n$ ; recall that a left
$R(G)$-module  
$B$ is called {\em nilpotent of class $\le k$}\/ if $I_R^{k}(G) \cdot B = 0$.  

\comment{
As usual, we shall write `$G$ -module' instead of
`$\Z(G)$-module'.
}

Now consider the map
  \[ \map{p_{n,R}}{G}{P_{n,R}(G,B)} \:,\quad p_{n,R}(a) = \rho(a-1) \:.\] 
In the case where $R=\Z$ we omit the subscript $R$. Recall that for a group homomorphism $f\mapco G \to H$ and an $R(H)$-module $M$, an $f$-{\em derivation}\/ from $G$ to $M$ is a map $d\mapco G\to M$ such that $d(ab)= f(a)d(b)+d(a)$ for $a,b\in G$.

\comment{
Note that for
$a,b\in G$,
\BE\label{wndpn} d_{p_{n,R}}(a,b) = w_n (p_{n-1,R}(a) \otimes p_{n-1,R}(b))\:.\EE
}

\begin{prop}\label{univpropspn} The maps $p_{n,R}$ and $p_n$ have the following universal properties:\V

(i) \map{p_{n,R}}{G}{P_{n,R}(G,B)} is a universal  $(G \auf G/B)$-derivation  from $G$  into  nilpotent
$R(G/B)$-modules of class $\le n$.\V

(ii) \map{p_{n}}{G}{P_{n}(G,B)} is a universal polynomial map of degree $\le n$ relative $B$ from $G$ into abelian groups.\V

\end{prop}

\N{\bf Proof\,:} (i) follows from the well known fact that the map $G\to I_{n,R}(G)$, $a\mapsto a-1$, is a universal  derivation 
from $G$ into arbitrary 
$G$-modules, see \cite{HiSt} VI.5. To prove (ii) we first show that $p_n$ is polynomial  of degree $\le n$ relative $B$. The linear
extension $\overline{p_n}\,\colon \Z(G) \to P_{n}(G,B)$ of $p_n$ satisfies $\overline{p_n}(a-1) = \rho(a-1)$ for $a\in G$, so $
\overline{p_n}(x) = \rho(x)$ for all $x\in I(G)$ by linearity. In particular, $\overline{p_n}(1+I(B)I(G)+I^{n+1}(G)) =
\rho(I(B)I(G)+I^{n+1}(G))= 0$, whence the assertion. To prove the universal property, let $f\,\colon G\to A$ be any polynomial map
of degree $\le n$ relative $B$. Then $\bar{f}\,\colon \Z(G)\to A$ factors through a \Z-linear map, also denoted by $\bar{f}$,
$\Z(G)/(I(B)I(G)+I^{n+1}(G)) \to A$. Denoting the restriction of $\bar{f}$ to $P_{n}(G,B)$ again by $\bar{f}$ we have
 $\bar{f}p_n = f$ as $\bar{f}(1) = 0$. Futhermore, $\bar{f}$ is the unique \Z-linear map with this  property since the elements
$p_n(a)$, $a\in G$, generate $P_{n}(G,B)$ as a \Z-module.\hfill$\Box$\V

Property (i) implies a  canonical isomorphism $\:P_{n,R}(G,B) 
\:\cong\: R(G/B) \otimes_{\Z(G)} P_{n}(G) \:$ of left $R(G/B)$-modules. 

Let $f\,\colon G\to A$ be a polynomial map of degree $\le n$ relative $B$ and $\bar{f}\,\colon\,
P_n(G,B)$ $ \to A$ the canonical induced \Z-linear map according to property (ii). Define the homomorphism
  \BE\label{wndef} w_f=\bar{f}\mu_n \,\colon\, P_{n-1}(G/B) \otimes_{\raisebox{-1mm}{$\sst \Z(G)$}} P_{n-1}(G,B) \:\lra \:
A\:.\EE
Then by \ref{polind} (1) we have the following commutative diagram of factorizations induced by $f$.

\BE\label{wfdia}
\begin{matrix}
G & \maprtop{f} & A & \mapltop{d_f} & G\times G \cr
\| & & \mapup{\bar{f}} & \nwarrow{\sst w_f}     & \mapdown{p_{n-1} \times p_{n-1}} \cr
G & \maprtop{p_n} & P_n(G,B) & \mapltop{\mu_n} &  P_{n-1}(G/B) \otimes_{\raisebox{-1mm}{$\sst \Z(G)$}} P_{n-1}(G,B) 
\end{matrix}\EE
\N Note that $w_{p_n} = \mu_n$.\V

Finally, suitable group homomorphisms induce ring
homomorphisms on the constructions introduced above in the obvious way.\V

From \REF{Iquot} and the inclusion $I_R(\gamma_n(G)) \subset I_R^n(G)$ we deduce the  natural isomorphisms of $R$-algebras
    \BE\label{Pn=Pn/gn} P_{n,R} \Big(G/\gamma_{n+1}(G)\,, B \gamma_{n+1}(G)/\gamma_{n+1}(G) \Big) \:\mapltop{\cong}\: P_{n,R}(G,B) \:\stackrel{\cong}{\lra}  P_{n,R}(G,B\gamma_n(G))
\EE


\nc{\ul}[1]{\underline{#1}}

Now let 
  \BE\label{Erw}  \ul{G}\colon\, B\Inj{i} G \Sur{\pi} Q \EE

\N be a group extension with abelian kernel
$B$. We will frequently identify $B$ with $i(B)$ and suppress $i$ from the notation. It is an elementary fact that the sequence

\BE\label{Pnrelsequ} R \ot_{\Z} B\:\maprtop{p_{n,R}\, i}\: P_{n,R}(G,B)
\Sur{P_{n,R}(\pi)} P_{n,R}(Q) \auf 0  \EE

\N is an exact sequence of $R(Q)$-linear homomorphisms, where the $Q$-action on $R \ot_{\Z}
 B$ is  
given by $R$-linear extension of the $Q$-action on $B$ induced by conjugation in $G$.
Moreover, the map $p_{n,R}\, i$ here denotes the $R$-linear extension of the map
$p_{n,R}\, i$  defined above.\V

\N Now consider the case where $B$ is {\em central}\/ in $G$. Then $I_R(B)I_R(G) = I_R(G)I_R(B)$, whence the map $\mu_n$ in
\REF{mundef} factors through another $R$-linear map, also denoted by $\mu_n$,
  \BE\label{muncentral}  \mu_n\,\colon\,  P_{n-1,R}(G/B) \ot_{\raisebox{-1mm}{$\sst R(G/B)$}}  P_{n-1,R}(G/B) \to P_{n,R}(G,B) \EE
 such that $\mu_n(p_{n-1,R}(a) \ot p_{n-1,R}(b)) = p_{n,R}(a)p_{n,R}(b)$, $a,b\in G$. Consequently, if \map{f}{G}{A} is a
polynomial map of degree $\le n$ relative $B$ then $w_f$ in \REF{wndef} factors through another \Z-linear map, also denoted by
$w_f$,
\BE\label{wfcentral}  w_f\,\colon\, P_{n-1}(G/B) \ot_{\raisebox{-1mm}{$\sst \Z(G/B)$}}  P_{n-1}(G/B) \to A \EE
such that $w_f(p_{n-1}(\bar{a}) \ot p_{n-1}(\bar{b})) = d_f(a,b)$, $a,b\in G$.

Now consider the case we are mainly interested in here, that is  $n=2$ and $B$ central in $G$. Noting that $P_1(G/B)$ is a trivial
$G/B$-module and using the canonical identifications
$P_1(G/B) \:\cong\: (G/B)\ab \:\cong\: G/BG\st$ we see that here $w_f$ is equivalent to the following linear map, also denoted by
$w_f$,
\BE\label{wf2}  w_f\,\colon\, \GBoGB \hspace{2mm} \stackrel{\cong}{\lra} \hspace{2mm} 
P_{1}(G/B) \ot_{\raisebox{-1mm}{$\sst \Z(G/B)$}}  P_{1}(G/B) \to A \:.\EE
This map satisfies $w_f(\bar{a}\ot \bar{b}) 
= d_f(a,b)$ for $a,b\in G$, whence it coincides with the map $w_f$ defined in Proposition
\ref{qumapsprops}(d) so that there is essentially no ambiguity in our notation.\V

Now we are ready for comparing the two universal constructions for quadratic and degree 2 polynomial maps, respectively.

\begin{prop}\label{Q=P2} Let $B$ be a central subgroup of a group $G$. Then there is a natural isomorphism of abelian groups
$\alpha\,\colon Q(G,B)\ab \stackrel{\cong}{\lra} P_2(G,B)$ such that $\alpha \rond ab \rond q = p_2$ and $\alpha \rond ab \rond w_q
= \mu_2$.
\end{prop}

\N{\bf Proof :} By \ref{quad=pol} the map \map{p_2}{G}{P_2(G,B)} is quadratic relative $B$ and the map
$ab\rond q \,\colon G \to Q(G,B)\ab$ is polynomial of degree $\le 2$. Hence $\alpha$ and its inverse are induced by the  
universal properties of $q$ and $p_2$, resp. So the equation $\alpha \rond ab \rond q = p_2$ holds by construction. Then the second
one follows from the definition of $w_q$ and from the identity $w_{p_n} = \mu_n$.\hfill $\Box$\V

Let $G$ be a group and $B$ a central subgroup of $G$. Then we have natural homomorphisms of abelian groups 
\BE \label{c2l2def} BG\st/\gamma_3(G) \hspace{2mm}\mapltop{c_2}\hspace{2mm} G/BG\st \sm G/BG\st \hspace{2mm}\Inj{l_2}
\hspace{2mm}\GBoGB \EE 
defined by $c_2(\bar{a} \sm \bar{b}) =  {[a,b]}\gamma_3(G)$ and $l_2(\bar{a} \sm \bar{b}) = \bar{a} \ot
\bar{b} - \bar{b} \ot
\bar{a}$.\V

\comment{
They are related by the equation
 \BE\label{c2l2} p_2ic_2 = \mu_2 l_2\EE
where  $i$
denotes the natural inclusion of $BG\st/\gamma_3(G)$ into $G/\gamma_3(G)$; this follows from the relation
 \BE\label{p2comm} p_2([a,b]) =  [p_2(a)\,,p_2(b)] = p_2(a)p_2(b) - p_2(b)p_2(a) \EE
for $a,b\in G$ which is immediate from equation \REF{comm-1}.\V
}

\begin{satz}\label{P2} Let $G$ be a group and $B$ a central subgroup of $G$. Then the following natural sequences of abelian groups
are exact.
\BE\label{GoGP2Gab} 0 \hspace{2mm} \lra \hspace{2mm} \Ker{c_2} \hspace{2mm} \maprtop{l_2} \hspace{2mm}  \GBoGB \hspace{2mm}
\maprtop{\mu_2} \hspace{2mm}  P_2(G,B)  \hspace{2mm} \maprtop{\rho_1} \hspace{2mm} G\ab \hspace{2mm} \lra \hspace{2mm} 1 \EE


\[ \makebox[0mm]{$0 \hspace{1mm} \to \hspace{1mm} G/BG\st \sm G/BG\st \hspace{2mm} \maprtop{(c_2,-l_2)^t}
\hspace{2mm} BG\st/\gamma_3(G) \hspace{2mm}\oplus\hspace{2mm} (\GBoGB\,)
\hspace{2mm}
\maprtop{(p_2i,\mu_2)} \hspace{2mm}  P_2(G,B)  \hspace{2mm} \maprtop{\rho_2} \hspace{2mm} G/BG\st \hspace{1mm} \to \hspace{1mm} 1$}
\]
\BE\label{pushP2GB} \EE
\BE\label{BP2P2} 0 \hspace{2mm} \lra \hspace{1mm} B\gamma_3(G)/\gamma_3(G) \hspace{2mm} \maprtop{p_2i} \hspace{2mm} P_2(G,B)
\hspace{2mm}
\maprtop{P_2(\pi)}
\hspace{2mm} P_2(G/B) \hspace{2mm} \lra \hspace{2mm} 0 \EE\vspace{2mm}

\N where $\rho_1p_2(a) =aG\st$ and $\rho_2p_2(a) =aBG\st$ for $a\in G$, $j\,\colon\, I^2(G)/I^3(G) \hra I(G)/I^3(G)= P_2(G)$, $i$
denotes the  injection of $BG\st/\gamma_3(G)$ or  $B\gamma_3(G)/\gamma_3(G)$ into $G/\gamma_3(G)$, and $\pi\,\colon\, (G,B)
\to (G/B,\{\bar{1}\})$ is the natural projection.\end{satz}\vspace{2mm}

\begin{bems} \rm Taking $B=\{1\}$ in \REF{GoGP2Gab} we rediscover the natural isomorphism $I^2(G)/I^3(G) \hcong {\rm U}_2{\rm L}(G)$ in \cite{Ba-GrII}. Moreover, centrality of $B$ is a crucial hypothesis in Theorem \ref{P2} as in general \REF{BP2P2} has to be replaced by the natural exact sequence 
 \[ {\rm Tor}_1^{\Z}(G/BG', G/BG') \mr{[\,,\,] \tau} B\gamma_3(G)/B'\gamma_3(G) \hspace{2mm} \maprtop{p_2i} \hspace{2mm} P_2(G,B)
\hspace{2mm}
\maprtop{P_2(\pi)}
\hspace{2mm} P_2(G/B) \hspace{2mm} \lra \hspace{2mm} 0\]
where the map $[\,,\,] \tau$ sends a typical generator $\langle aBG',k,bBG' \rangle$ with $a,b\in G$, $k\in \Z$ such that $a^k,b^k\in BG'$ (see \cite[V.6]{ML}), to the element $[a,b^k]B'\gamma_3(G)$ which is nontrivial in general, see Theorem 2.6 and Example 2.4 in \cite{D3F2}. 
\end{bems}

\N{\bf Proof of Theorem \ref{P2}\,:} By the right hand isomorphism in \REF{Pn=Pn/gn} we may assume that $\gamma_3(G)=1$. Now recall the central
extension 
\[ 0 \:\lra\: \GBoGB \hspace{2mm} \maprtop{w_q} \hspace{2mm}  Q(G,B) \hspace{2mm} \maprtop{\widehat{id}} \hspace{2mm} G \:\lra\:
1 \]
from \ref{Q(G)}. Putting $P=\widehat{id}^{-1}G\st$ and writing $[-,-]$ for the respective commutator maps we have the following
commutative diagram with exact rows.

\nc{\incdown}[1]{@\cdvgeneric>\kern-1pt\scriptscriptstyle\cap\kern2.5pt> 
\arrowvert>\downarrow>  > inc >} 

\[ \begin{CD}
\begin{matrix}
 & & Q(G,B) \times Q(G,B) & \maprtop{\widehat{id}\times \widehat{id}} & G\times G\ruleu\ruled \cr
 & & \mapdown{[-,-]} & \mapsw{c} & \mapdown{[-,-]} \cr
\GBoGB & \Inj{w_q} & P & \Sur{\widehat{id}} & G\st \ruleu\cr
\ruleu\| & & \incdown{inc}   \incdown{inc} \cr
\GBoGB & \Inj{w_q} & Q(G,B) & \Sur{\widehat{id}} & G \ruled\ruleu \cr
\end{matrix}\end{CD}\]
The factorisation through $c$ of the commutator map of $Q(G,B)$ exists as \Imm{w_q} is central in $Q(G,B)$; it satisfies 
\BE\label{c=} c(a,b) = [q(a),q(b)] = q[a,b] - w_q(\bar{a}\ot \bar{b} - \bar{b}\ot \bar{a}) \EE
by \REF{f[a,b]}.  As 
$[-,-]\,\colon\,G\times G \to G\st$ is bilinear and $q$ is linear on $G\st$  by Proposition \ref{qumapsprops}(c), $c$ is bilinear; it annihilates
$BG\st$ in both variables since $B$ is central and $G\st$ is abelian. Now $P$ is abelian being a split central extension of the
abelian group
$G\st$ by construction of $Q(G,B)$. Hence $c$ gives rise to a homomorphism $\bar{c}\,\colon\,G/BG\st \sm G/BG\st \to P$ such that
$\bar{c}(\bar{a} \sm \bar{b}) = [q(a),q(b)]$ and $Q(G,B)\st = \langle \Imm{c} \rangle = \Imm{\bar{c}}$. Consider the following
commutative diagram with exact rows.
\[
\begin{matrix}
\GBoGB & \Inj{w_q} & P & \Sur{\widehat{id}} & G\st \ruled\ruleu \cr
\mapup{-\,l_2} & & \mapup{\bar{c}} & & \| \cr
\Ker{c_2} & \Inj{inc} & G/BG\st \sm G/BG\st & \Sur{c_2} & G\st\ruled
\end{matrix}\]
The right hand square clearly commutes, and by \REF{c=}, $\bar{c} = qc_2 - w_ql_2$, whence the left hand square commutes, too.
It follows that $Q(G,B)\ab = {\rm coker}(inc \rond \bar{c})$ fits into the exact sequence
 \[ 0 \:\to (\GBoGB) \Big/ l_2\Ker{c_2} \hspace{2mm} \maprtop{ab\rond w_q} \hspace{2mm} Q(G,B)\ab \hspace{2mm}\lra \hspace{2mm} G\ab
\: \to 1\] 
which becomes sequence \REF{GoGP2Gab} under the identification $Q(G,B)\ab \:\cong\:P_2(G,B)$, see \ref{Q=P2}. 
As to sequence \REF{pushP2GB}, first note that it is exact in $P_2(G,B)$ since the composite
isomorphism ${\rm coker}(\mu_2) \:\cong\: I(G)/I^2(G)  \:\cong\: G/G\st$ takes $\overline{p_2i(b)}$ to $bG\st$ for $b\in B$.
Furthermore, we have $(p_2i,\mu_2) (c_2,-l_2)^t = 0$ by the relation $ p_2ic_2 = \mu_2 l_2$
which follows from the identity
 \BE\label{p2comm} p_2([a,b]) =  [p_2(a)\,,p_2(b)] = p_2(a)p_2(b) - p_2(b)p_2(a) \EE
for $a,b\in G$; this is immediate from equation \REF{comm-1}.
 The
map $(c_2,-l_2)^t$ is injective as $l_2$ is; it remains to show that the map $\overline{(p_2i,\mu_2)}\,\colon\, \Pi \colon= {\rm
coker}(c_2,-l_2)^t \to P_2(G,B)$ induced by $(p_2i,\mu_2)$ is injective. Consider the  commutative diagram \REF{P2bewdia} below
where $\phi$ is the natural projection and $i_1,i_2$ are the natural inclusions into $BG\st \:\oplus\: (\GBoGB\,)$ followed by the
natural projection to $\Pi$.

  \BE\label{P2bewdia}
\begin{CD}
\begin{matrix}
G/BG\st \sm G/BG\st & \maprtop{c_2} & BG\st & \Sur{\phi} & BG\st/G\st \ruleu\ruled\cr
\mapdown{l_2} & & \mapdown{i_1} & & \| \cr
\GBoGB & \maprtop{i_2} & \Pi & \Sur{\overline{(\phi,0)}} & BG\st/G\st \ruled \cr
\surdown{} & & \| & & \| \cr
\GBoGB\Big/ l_2\Ker{c_2}  & \Inj{\overline{i_2}} & \Pi & \Sur{\overline{(\phi,0)}} & BG\st/G\st \cr
\| & & \mapdown{\overline{(p_2i,\mu_2)}} & & \incdown{inc} \cr
\GBoGB\Big/ l_2\Ker{c_2} & \Inj{\overline{\mu_2}} & P_2(G,B) & \Sur{\rho_1} & G\ab \rule[-5mm]{0mm}{0mm}
\end{matrix}
\end{CD}
\EE
The rows are exact: for the bottom row this follows from sequence \REF{GoGP2Gab}, for the second and third row from the fact that the upper left hand square is a cocartesian square
of abelian groups (or by easy direct arguments).
Injectivity of $\overline{(p_2i,\mu_2)}$ is now immediate. 

Finally, to prove exactness of sequence \REF{BP2P2} it suffices to check injectivity of $p_2i$, see \REF{Pnrelsequ}. But in the
above diagram $i_1$ is injective as $l_2$ is, so $p_2i = \overline{(p_2i,\mu_2)} \rond i_1$ is injective, too.\hfill$\Box$\V

\tw